\documentclass{elsart}
\usepackage{amssymb}
\usepackage{pstricks, pst-plot, pst-node, pst-tree}
\usepackage{array}
\usepackage{fancyhdr}
\usepackage{longtable}
\usepackage{tabularx}
\usepackage{rotating}

\usepackage{makeidx}
\usepackage{amsmath}
\usepackage{theorem}
\usepackage[mathscr]{eucal}
\usepackage{enumerate}
\usepackage{float}
\usepackage{ifthen}
\newif\ifPDF
\ifx\pdfoutput\undefined
  \PDFfalse
\else
  \PDFtrue
\fi
\ifPDF
\usepackage[pdftex]{epsfig}
\usepackage[pdftex]{graphicx}
\else
\usepackage[dvips]{epsfig}
\usepackage[dvips]{graphicx}
\fi
\journal{Journal of Computational and Applied Mathematics}
\volume{214}
\pubyear{2008}
\newcounter{mylastpage}
\issue{1}
\usepackage{hyperref}
\makeatletter
\def\ps@copyright{%
 \let\@oddhead\@empty
 \let\@evenhead\@empty
 \def\@oddfoot{\scriptsize\fontsize{10}{13}\selectfont\slshape\hskip-7em
   Published in \@journal\ \@volume\ (\the\@pubyear)\ no.\ \@issue, pp.\ \ESpagenumber{firstpage}--\ESpagenumber{mylastpage},\newline
\href{http://dx.doi.org/10.1016/j.cam.2007.02.040}{doi: 10.1016/j.cam.2007.02.040}}%
 \let\@evenfoot\@oddfoot}
\makeatother
\theoremstyle{plain}
\newtheorem{Def}{Definition}[section]

\newtheorem{The}[Def]{Theorem}

\newtheorem{Bem}[Def]{Remark}

\begin{document}

\begin{frontmatter}
\title{Continuous Weak Approximation for
Stochastic Differential Equations}
\author{Kristian Debrabant} and
\ead{debrabant@mathematik.tu-darmstadt.de}
\author{Andreas R\"{o}{\ss}ler}
\ead{roessler@mathematik.tu-darmstadt.de}
\address{Technische Universit\"{a}t Darmstadt, Fachbereich Mathematik, Schlo{\ss}gartenstr.7,
D-64289 Darmstadt, Germany}

\begin{abstract}
A convergence theorem for the continuous weak approximation of the
solution of stochastic differential equations by general one step
methods is proved, which is an extension of a theorem due to
Milstein. As an application, uniform second order conditions for a
class of continuous stochastic Runge--Kutta methods containing the
continuous extension of the second order stochastic Runge--Kutta
scheme due to Platen are derived. Further, some coefficients for
optimal continuous schemes applicable to It\^{o} stochastic
differential equations with respect to a multi--dimensional Wiener
process are presented.
\end{abstract}

\end{frontmatter}

\section{Introduction} \label{Introduction}
Since in recent years the application of stochastic differential
equations (SDEs) has increased rapidly, there is now also an
increasing demand for numerical methods. Many numerical schemes
have been proposed in literature, see e.g., Kloeden and
Platen~\cite{KP99} or Milstein and Tretyakov~\cite{Mil04} and the
references therein. The present paper deals with conditions on the
local error of continuous one step methods for the approximation
of the solution $X(t)$ of stochastic differential equations such
that global convergence in the weak sense is assured. The
continuous methods under consideration are also called dense
output formulas~\cite{HNW93}. Continuous methods are applied
whenever the approximation $Y(t)$ has to be determined at
prescribed dense time points which would require the step size to
be very small. Also if a graphical output of the approximate
solution is needed, continuous methods may be applied as well.
Further, for the weak approximation of the solution of stochastic
differential delay equations with variable step size, a global
approximation to the solution is needed. In the deterministic case
it is known~\cite{HNW93} that continuous methods may be superior
over many other interpolation methods. Therefore, the application
of continuous methods for stochastic differential delay equations
may be very promising for future research (see also
\cite{BucShard05}). \\ \\
As an example, the continuous extension of a certain class of
stochastic Runge-Kutta methods is given. Some stochastic
Runge--Kutta (SRK) methods for strong ap\-pro\-xi\-ma\-tion have
been introduced by Burrage and Burrage~\cite{BuBu00a,BuBu96} and
Newton~\cite{New91}. For the weak approximation, SRK methods have
been developed as well, see e.g., Komori and Mitsui~\cite{KoMi95},
Komori, Mitsui and Sugiura~\cite{KoMiSu97},
R\"{o}{\ss}ler~\cite{Roe04,Roe04d,Roe04e,Roe03} or Tocino and
Ardanuy~\cite{ToA02} and Tocino and Vigo-Aguiar~\cite{ToVA02}.
\\ \\
The main advantage of continuous SRK methods is the cheap
numerical approximation of $E(f(X(t_0 + \theta h)))$ for the whole
integration interval $0 \leq \theta \leq 1$ beside the
approximation of $E(f(X(t_0+h)))$. Here, cheap means without
additional evaluations of drift and diffusion and without the
additional simulation of random variables. \\ \\
The paper is organized as follows: We first prove a convergence
theorem for the uniform continuous weak approximation by one step
methods in Section~\ref{St-SubSection-Local-Global-Approx}, which
is a modified version of a theorem due to Milstein~\cite{Mil95}.
Based on the convergence theorem, we extend a class of SRK methods
of weak order one and two to continuous stochastic Runge-Kutta
(CSRK) methods. Here, the considered class of SRK methods contains
the well known weak order two SRK scheme due to
Platen~\cite{KP99}. Further, we give order conditions for the
coefficients of this class of CSRK methods in
Section~\ref{SRK-methods-Ito-SDE-systems-Heading1}. Taking into
account additional order conditions in order to minimize the
truncation error of the approximation yields optimal CSRK schemes.
Finally, the performance of one optimal CSRK scheme is analyzed by
a numerical example in Section~\ref{Sec:Numerical-Example}.
\section{Local and Global Weak Convergence}
\label{St-SubSection-Local-Global-Approx}
Let $(X(t))_{t \in I}$ denote the solution process of a
$d$-dimensional It\^{o} stochastic differential equation
\begin{equation} \label{St-lg-sde-ito-1}
    {\mathrm{d}} X(t) = a(t,X(t)) \, {\mathrm{d}}t + b(t,X(t)) \,
    {\mathrm{d}}W(t), \qquad X(t_0) = X_{0},
\end{equation}
with a driving $m$-dimensional Wiener process $(W(t))_{t \geq 0}$
and with $I=[t_0,T]$.
It is always assumed that the Borel-measurable coefficients $a : I
\times \mathbb{R}^d \rightarrow \mathbb{R}^d$ and $b : I \times
\mathbb{R}^d \rightarrow \mathbb{R}^{d \times m}$ satisfy a
Lipschitz and a linear growth condition
\begin{alignat}{2} \label{lipschitz}
    &\| a(t,x)-a(t,y) \| + \| b(t,x)-b(t,y) \| \leq C \, \| x-y
    \|, \\
    \label{linear-growth}
    &\| a(t,x) \|^2 + \| b(t,x) \|^2 \leq C^2 (1+\|x\|^2) ,
\end{alignat}
for every $t \in I$, $x,y \in \mathbb{R}^d$ and a constant $C>0$
such that the Existence and Uniqueness Theorem~4.5.3~\cite{KP99}
applies. Here, $\| \cdot \|$ denotes a vector norm, e.g., the
Euclidian norm, or the corresponding matrix norm. \\ \\
Let a discretization $I_h = \{t_0, t_1, \ldots, t_N\}$ with $t_0 <
t_1 < \ldots < t_N =T$ of the time interval $I=[t_0,T]$ with step
sizes $h_n = t_{n+1}-t_n$ for $n=0,1, \ldots, N-1$ be given.
Further, define $h = \max_{0 \leq n < N} h_n$ as the maximum step
size.
Let $X^{t_0,X_{0}}$ denote the solution of the stochastic
differential equation (\ref{St-lg-sde-ito-1}) in order to
emphasize the initial condition. If the coefficients of the
It\^{o} SDE (\ref{St-lg-sde-ito-1}) are globally Lipschitz
continuous, then there always exists a version of $X^{s,x}(t)$
which is continuous in $(s,x,t)$ such that $X^{t_0,X_{0}}(t) =
X^{s,X^{t_0,X_{0}}(s)}(t)$ $P\text{-a.s.}$ holds for $s,t \in I$
with $s \leq t$ (see \cite{Kun84}).
For simplicity of notation, in this section it is supposed that
$I_h$ denotes an equidistant discretization, i.e.\ $h =
(t_N-t_0)/N$. Next, we consider the one-step approximation
\begin{equation} \label{St-lg-approxi}
    Y^{t,x}(t+h) = A(t,x,h; \xi),
\end{equation}
where $\xi$ is a vector of random variables, with moments of
sufficiently high order, and $A$ is a vector function of dimension
$d$. We write $Y_{n+1} = Y^{t_n,Y_n}(t_{n+1})$ and we construct
the sequence
\begin{equation} \label{St-lg-method}
    \begin{split}
    Y_0 &= X_0, \\
    Y_{n+1} &= A(t_n, Y_n, h;
    \xi_n), \quad n=0,1, \ldots,N-1,
    \end{split}
\end{equation}
where $\xi_0$ is independent of $Y_0$, while $\xi_n$ for $n \geq
1$ is independent of $Y_0, \ldots, Y_n$ and $\xi_0, \ldots,
\xi_{n-1}$. Then $Y = \{Y^{s,x}(t) : x \in \mathbb{R}, s \leq t
\text{ with } s,t \in I_h\}$ is a (non-homogeneous) discrete time
Markov chain depending on $h$ such that $Y^{t_0,X_0}(t) =
Y^{s,Y^{t_0,X_0}(s)}(t)$ $P\text{-a.s.}$ holds for $s,t \in I_h$
with $s \leq t$. \\ \\
In the following, let $C_P^l(\mathbb{R}^d, \mathbb{R})$ denote the
space of all $g \in C^l(\mathbb{R}^d,\mathbb{R})$ with polynomial
growth, i.e.\ there exists a constant $C > 0$ and $r \in
\mathbb{N}$, such that $| \partial_x^i g(x) | \leq C
(1+\|x\|^{2r})$ holds for all $x \in \mathbb{R}^d$ and any partial
derivative of order $i \leq l$. Further, let  $g \in C_P^{k,l}(I
\times \mathbb{R}^d, \mathbb{R})$ if $g(\cdot,x) \in
C^{k}(I,\mathbb{R})$ and $g(t,\cdot) \in C_P^l(\mathbb{R}^d,
\mathbb{R})$ for all $t \in I$ and $x \in \mathbb{R}^d$
\cite{KP99}. \\ \\
Since we are interested in obtaining a continuous global
approximation converging in the weak sense with some desired order
$p$, we give an extension of the convergence theorem due to
Milstein~\cite{Mil95,Mil04} which specifies the relationship
between the local and the global approximation order.

\begin{The} \label{St-lg-theorem-main}
    Suppose the following conditions hold:
    \begin{enumerate}[(i)]
        \item \label{St-lg-cond1} The coefficients $a^i$
        and $b^{i,j}$ are continuous,
        satisfy a Lipschitz condition (\ref{lipschitz}) and belong to
        $C_P^{p+1,2(p+1)}(I \times \mathbb{R}^d, \mathbb{R})$ with respect to $x$ for $i=1, \ldots,
        d$, $j=1, \ldots, m$.
        \item \label{St-lg-cond2} For sufficiently large $r$ (specified in the proof)
        the moments $E(\|Y_n\|^{2r})$ exist and are uniformly bounded with respect
        to $N$ and $n=0,1, \ldots, N$.
        \item \label{St-lg-cond3} Assume that for all $f \in
        C_P^{2(p+1)}(\mathbb{R}^d, \mathbb{R})$ there exists a $K \in
        C_P^0(\mathbb{R}^d, \mathbb{R})$ such that the following {\emph{local error
        estimations}}
        \begin{eqnarray}
            |E(f(X^{t,x}(t+h))) - E(f(Y^{t,x}(t+h)))| &\leq& K(x) \, h^{p+1}
            \label{St-lg-local-error-estimation:a}\\
            |E(f(X^{t,x}(t+\theta h))) - E(f(Y^{t,x}(t+\theta h)))| &\leq& K(x) \,
            h^{p} \label{St-lg-local-error-estimation:b}
        \end{eqnarray}
        are valid for $x \in \mathbb{R}^d$, $t, t+h \in [t_0,T]$ and $\theta \in [0,1]$.
    \end{enumerate}
    Then for all $t \in [t_0,T]$ the following
    {\emph{global error estimation}}
    \begin{equation} \label{St-lg-global-error-estimation}
        | E(f(X^{t_0,X_0}(t))) - E(f(Y^{t_0,X_0}(t))) | \leq C \, h^p
    \end{equation}
    holds for all $f \in C_P^{2(p+1)}(\mathbb{R}^d, \mathbb{R})$, where $C$ is a
    constant, i.e.\ the method (\ref{St-lg-method}) has a uniform order of
    accuracy $p$ in the sense of weak approximation.
\end{The}

\begin{Bem}
    In contrast to the original theorem now the order of
    convergence specified in
    equation~(\ref{St-lg-global-error-estimation}) is not only
    valid in the discretization times $t \in I_h$. Provided that
    the additional
    condition~(\ref{St-lg-local-error-estimation:b})
    is fulfilled, the global order of convergence
    (\ref{St-lg-global-error-estimation})
    holds also uniformly for all $t \in [t_0,T]$.
\end{Bem}

{\textbf{Proof of Theorem~\ref{St-lg-theorem-main}.}} The proof
extends the ideas of the proof of Theorem~9.1 in \cite{Mil95}.
However, now all time points $t \in [t_0,T]$ have to be considered
instead of only $t \in I_h$ like in \cite{Mil95} and an additional
estimate is necessary in order to prove the uniform order of
convergence on the whole time interval. Therefore, we consider the
function
\begin{equation*}
    u(s,x) = E(f(X(t_{k}+\theta h)) | X(s)=x )
\end{equation*}
for $s, t_k+ \theta h \in I$, $x \in \mathbb{R}^d$ and $t_{k} \in
I_h$ with $s \leq t_{k}$ and for $\theta \in [0,1]$. Due to
condition~(\ref{St-lg-cond1}) $u(s, \cdot)$ has partial
derivatives of order up to $2(p+1)$, inclusively, and belongs to
$C_P^{2(p+1)}(\mathbb{R}^d, \mathbb{R})$ for each $s \in [t_0,
t_k]$ \cite{KP99,Mil95}. Therefore, $u$ satisfies
\begin{equation*}
    \begin{split}
    |E(u(s ,X^{t,x}(t+h))) - E(u(s ,Y^{t,x}(t+h)))| &\leq K_u(x) \, h^{p+1}
    \end{split}
\end{equation*}
uniformly w.r.t.\ $s \in [t_0, t_{k}]$ for some $K_u \in
C_P^{0}(\mathbb{R}^d, \mathbb{R})$. Since $Y_0 = X_0$, we have
\begin{equation} \label{St-lg-proof-eq2}
    \begin{split}
    E(f(X^{t_0,X_0}(t_{k}+\theta h))) = & \, \, E(f(X^{t_1,X^{t_0,Y_0}(t_1)}(t_{k}+\theta h))) -
    E(f(X^{t_1,Y_1}(t_{k}+\theta h))) \\
    & + E(f(X^{t_1,Y_1}(t_{k}+\theta h))).
    \end{split}
\end{equation}
Furthermore, since $X^{t_1,Y_1}(t_{k}+\theta h) = X^{t_2,
X^{t_1,Y_1}(t_2)}(t_{k}+\theta h)$, we have
\begin{equation} \label{St-lg-proof-eq3}
    \begin{split}
    E(f(X^{t_1,Y_1}(t_{k}+\theta h))) = & \, \,
    E(f(X^{t_2,X^{t_1,Y_1}(t_2)}(t_{k}+\theta h))) -
    E(f(X^{t_2,Y_2}(t_{k}+\theta h))) \\
    & + E(f(X^{t_2,Y_2}(t_{k}+\theta h))).
    \end{split}
\end{equation}
By (\ref{St-lg-proof-eq2}) and (\ref{St-lg-proof-eq3}) we get
\begin{equation*}
\begin{split}
    E(f(X^{t_0,X_0}&(t_{k}+\theta h))) \\
    = &
    E(f(X^{t_1,X^{t_0,Y_0}(t_1)}(t_{k}+\theta h))) -
    E(f(X^{t_1,Y_1}(t_{k}+\theta h))) \\
     & + E(f(X^{t_2,X^{t_1,Y_1}(t_2)}(t_{k}+\theta h))) -
    E(f(X^{t_2,Y_2}(t_{k}+\theta h))) \\
    & + E(f(X^{t_2,Y_2}(t_{k}+\theta h))).
\end{split}
\end{equation*}
The procedure continues to obtain
\begin{equation*}
    \begin{split}
    E(f(X^{t_0,X_0}&(t_{k}+\theta h))) \\
    = &\sum_{i=0}^{k-1} \left(
    E(f(X^{t_{i+1},X^{t_i,Y_i}(t_{i+1})}(t_{k}+\theta h))) -
    E(f(X^{t_{i+1},Y_{i+1}}(t_{k}+\theta h))) \right) \\
    &+ E(f(X^{t_k,Y_k}(t_{k}+\theta h))).
    \end{split}
\end{equation*}
Recalling that $Y^{t_0,X_0}(t_{k} + \theta h) = Y^{t_k,Y_k}(t_{k}
+ \theta h)$, this implies the identity
\begin{alignat}{2} \label{St-lg-proof-eq6}
    E(f(X^{t_0,X_0}(t_{k}+\theta h))) &- E(f(Y^{t_0,X_0}(t_{k}+\theta h))) \notag \\
    =&\sum_{i=0}^{k-1}  \left(
    E(E(f(X^{t_{i+1},X^{t_i,Y_i}(t_{i+1})}(t_{k}+\theta h)) \, | \,
    X^{t_i,Y_i}(t_{i+1}))) \right. \notag \\
    & \left. - E(E(f(X^{t_{i+1},Y^{t_i,Y_i}(t_{i+1})}(t_{k}+\theta h))
    \, | \, Y^{t_i,Y_i}(t_{i+1}))) \right) \notag \\
    & + E(f(X^{t_k,Y_k}(t_{k}+\theta h))) - E(f(Y^{t_k,Y_k}(t_{k}+\theta h))).
\end{alignat}
According to the definition of $u(t,x)$, the Jensen inequality and
(\ref{St-lg-proof-eq6}) imply
\begin{alignat}{2} \label{St-lg-proof-eq7}
    | E(f(X^{t_0,X_0}&(t_{k}+\theta h))) - E(f(Y^{t_0,X_0}(t_{k}+\theta h))) |
    \notag \\
    = & \Big| \sum_{i=0}^{k-1} \big( E(u(t_{i+1},
    X^{t_i,Y_i}(t_i+h))) - E(u(t_{i+1},Y^{t_i,Y_i}(t_i+h)))
    \big) \notag \\
    & + \big( E(f(X^{t_k,Y_k}(t_{k}+\theta h))) -
    E(f(Y^{t_k,Y_k}(t_{k}+\theta h))) \big) \Big| \notag \\
    \leq & \sum_{i=0}^{k-1} E \left( \big| E \big( u(t_{i+1},
    X^{t_i,Y_i}(t_i+h)) - u(t_{i+1}, Y^{t_i,Y_i}(t_i+h)) \, | \,
    Y_i \big) \big| \right) \notag \\
    & + E \left( \big| E\big(f(X^{t_k, Y_k}(t_{k}+\theta h)) -
    f(Y^{t_k,Y_k}(t_{k}+\theta h)) \, | \, Y_k \big) \big| \right).
\end{alignat}
Notice that the functions $u(t,\cdot)$ and $f$, which belong to
$C_P^{2(p+1)}(\mathbb{R}^d, \mathbb{R})$ and satisfy an inequality
of the form (\ref{St-lg-local-error-estimation:a}) and
(\ref{St-lg-local-error-estimation:b}), also satisfy along with it
a conditional version of such an inequality. Let $K \in
C_P^0(\mathbb{R}^d, \mathbb{R})$ such that for both $u(t,\cdot)$
and $f$ the inequalities (\ref{St-lg-local-error-estimation:a})
and (\ref{St-lg-local-error-estimation:b}) hold. Thus there exist
$r \in \mathbb{N}$ and $C > 0$ such that
\begin{equation} \label{St-lg-proof-eq8}
    | K(x) | \leq C (1+ \| x \|^{2r})
\end{equation}
holds for all $x \in \mathbb{R}^d$. Then (\ref{St-lg-proof-eq7})
together with (\ref{St-lg-proof-eq8}) imply for $0 < \theta h \leq
h$ the estimate
\begin{alignat}{2} \label{St-lg-proof-eq9}
    | E(f(&X^{t_0,X_0}(t_{k}+\theta h))) - E(f(Y^{t_0,X_0}(t_{k}+\theta h))) |
    \notag \\
    \leq & \sum_{i=0}^{k-1} E( | K(Y_i) \, h^{p+1} | ) + E( |
    K(Y_k) \, (\theta h)^{p} | ) \notag \\
    \leq & \sum_{i=0}^{k-1} C (1 + E( \| Y_i \|^{2r})) \,
    h^{p+1} + C (1 + E( \| Y_k \|^{2r})) \, (\theta h)^{p} \notag \\
    \leq & k \, C (1 + \max_{0 \leq i \leq k-1} E( \| Y_i \|^{2r})) \,
    h^{p+1} + C (1 + E( \| Y_k \|^{2r})) \, h^{p}.
\end{alignat}
Assuming that condition (\ref{St-lg-cond2}) holds precisely for
this $2r$ and applying finally $k = \tfrac{t_{k} - t_0}{h}$ in
(\ref{St-lg-proof-eq9}), (\ref{St-lg-global-error-estimation}) is
obtained. \hfill $\square$ \\ \\
In the following, we assume that the coefficients $a^i$ and
$b^{i,j}$ satisfy assumption~(\ref{St-lg-cond1}) of
Theorem~\ref{St-lg-theorem-main}. Further,
assumption~(\ref{St-lg-cond2}) of Theorem~\ref{St-lg-theorem-main}
is always fulfilled for the class of stochastic Runge--Kutta
methods considered in the present paper (see \cite{Roe04d,Roe03}
for details).
\section{Continuous Stochastic Runge--Kutta Methods}
\label{SRK-methods-Ito-SDE-systems-Heading1}
As an example, we consider the continuous extension of the class
of stochastic Runge--Kutta methods due to
R\"{o}{\ss}ler~\cite{Roe04d,Roe04e,Roe03} which contains the weak
order two Runge--Kutta type scheme due to Platen~\cite{KP99}. The
intention is to approximate the solution of the It\^{o}
SDE~(\ref{St-lg-sde-ito-1}) on the whole interval $I=[t_0,T]$ in
the weak sense. Therefore, we define the $d$-dimensional
approximation process $Y$ by an explicit continuous stochastic
Runge--Kutta method having $s$ stages with initial value $Y(t_0) =
X_0$ and
\begin{equation} \label{SRK-method-Ito-Wm-allg01}
    \begin{split}
    Y(t_n+\theta \, h_n) &= Y(t_n) + \sum_{i=1}^s
    \alpha_i(\theta) \,\, a(t_n+c_i^{(0)} h_n, H_i^{(0)}) \, h_n \\
    & +
    \sum_{i=1}^s \sum_{k=1}^m
    \left( \beta_i^{(1)}(\theta) \, \hat{I}_{(k)} + \beta_i^{(2)}(\theta) \,
    \tfrac{\hat{I}_{(k,k)}}{\sqrt{h_n}} \right) b^{k}(t_n+c_i^{(1)} h_n, H_i^{(k)}) \\
    & + \sum_{i=1}^s
    \sum_{\substack{k,l=1 \\ k \neq l}}^m
    \left( \beta_i^{(3)}(\theta) \, \hat{I}_{(k)} + \beta_i^{(4)}(\theta) \,
    \tfrac{\hat{I}_{(k,l)}}{\sqrt{h_n}} \right)
    b^{k}(t_n+c_i^{(2)} h_n, \hat{H}_i^{(l)})
    \end{split}
\end{equation}
for $\theta \in [0,1]$ and $n=0,1, \ldots, N-1$ with stage values
\begin{alignat*}{5}
    H_i^{(0)} &=&\,\, Y(t_n) &+ \sum_{j=1}^{i-1} A_{ij}^{(0)}
    \, a(t_n+c_j^{(0)} h_n, H_j^{(0)}) \, h_n \\
    && &+ \sum_{j=1}^{i-1} \sum_{r=1}^m
    B_{ij}^{(0)} \, b^r(t_n+c_j^{(1)} h_n, H_j^{(r)}) \, \hat{I}_{(r)} \\
    H_i^{(k)} &=&\,\, Y(t_n) &+ \sum_{j=1}^{i-1} A_{ij}^{(1)}
    \, a(t_n+c_j^{(0)} h_n, H_j^{(0)}) \, h_n \\
    && &+ \sum_{j=1}^{i-1}
    B_{ij}^{(1)} \, b^k(t_n+c_j^{(1)} h_n, H_j^{(k)}) \,
    \sqrt{h_n} \\
    \hat{H}_i^{(k)} &=&\,\, Y(t_n) &+ \sum_{j=1}^{i-1} A_{ij}^{(2)}
    \, a(t_n+c_j^{(0)} h_n, H_j^{(0)}) \, h_n \\
    && &+ \sum_{j=1}^{i-1}
    B_{ij}^{(2)} \, b^k(t_n+c_j^{(1)} h_n, H_j^{(k)}) \, \sqrt{h_n}
\end{alignat*}
for $i=1, \ldots, s$ and $k=1, \ldots, m$. Here, the weights
$\alpha_i, \beta_i^{(r)} \in C([0,1],\mathbb{R})$ are some
continuous functions for $1 \leq i \leq s$. We denote
$\alpha(\theta) = (\alpha_1(\theta), \ldots, \alpha_s(\theta))^T$,
$\beta^{(r)}(\theta)=(\beta_1^{(r)}(\theta), \ldots,
\beta_s^{(r)}(\theta))^T$, $c^{(q)} \in \mathbb{R}^s$ and
$A^{(q)}, B^{(q)} \in \mathbb{R}^{s \times s}$ for $0 \leq q \leq
2$, $1 \leq r \leq 4$ and with $A_{ij}^{(q)} = B_{ij}^{(q)} =0$
for $j \geq i$, which are the vectors and matrices of coefficients
of the SRK method. We choose $c^{(q)}=A^{(q)} e$ for $0 \leq q
\leq 2$ with a vector $e=(1, \ldots, 1)^T$ \cite{Roe04d}. In the
following, the product of vectors is defined component-wise.
The coefficients of the CSRK method can be arranged in the
following Butcher tableau:

\renewcommand{\arraystretch}{1.8}
\begin{center}
\begin{tabular}{c|c|c|c}
    $c^{(0)}$ & ${A}^{(0)}$ & $B^{(0)}$ & \\
    \cline{1-4}
    $c^{(1)}$ & ${A}^{(1)}$ & $B^{(1)}$ & \\
    \cline{1-4}
    $c^{(2)}$ & ${A}^{(2)}$ & $B^{(2)}$ & \\
    \hline
    & $\alpha^T$ & ${\beta^{(1)}}^T$ & ${\beta^{(2)}}^T$ \\
    \cline{2-4}
    & & ${\beta^{(3)}}^T$ & ${\beta^{(4)}}^T$
\end{tabular}
\end{center}
\renewcommand{\arraystretch}{1.0}

The random variables are defined by $\hat{I}_{(r)} = \Delta
\hat{W}_n^{r}$ and $ \hat{I}_{(k,l)} = \frac{1}{2} ( \Delta
\hat{W}_n^{k} \, \Delta \hat{W}_n^{l} + V_n^{k,l} )$. Here, the
$\Delta \hat{W}_n^r$ are independently three-point distributed
random variables with $P(\Delta \hat{W}_n^r = \pm \sqrt{3 \, h_n}
) = \frac{1}{6}$ and $P(\Delta \hat{W}_n^r = 0 ) = \frac{2}{3}$
for $1 \leq r \leq m$. Further, the $V_n^{k,l}$ are independently
two-point distributed random variables with $P(V_n^{k,l} = \pm
h_n) = \frac{1}{2}$ for $k=1, \ldots, k-1$, $V_n^{k,k} = -h_n$ and
$V_n^{k,l} = -V_n^{l,k}$ for $l=k + 1, \ldots, m$ and $k=1,
\ldots, m$ \cite{KP99}. \\ \\
The algorithm works as follows: First, the random variables
$\hat{I}_{(k)}$ and $\hat{I}_{(k,l)}$ have to be simulated for
$1\leq k,l \leq m$ w.r.t.\ the actual step size $h_n$. Next, based
on the approximation $Y(t_n)$ and the random variables, the stage
values $H^{(0)}$, $H^{(k)}$ and $\hat{H}^{(k)}$ are calculated.
Then we can determine the continuous approximation $Y(t)$ for
arbitrary $t \in [t_n,t_{n+1}]$ by varying $\theta$ from 0 to 1 in
formula~(\ref{SRK-method-Ito-Wm-allg01}). Thus, only a very small
additional computational effort is needed for the calculation of
the values $Y(t)$ with $t \in I \setminus I_h$. This is the main
advantage in comparison to the application
of an SRK method with very small step sizes. \\ \\
Using the multi--colored rooted tree analysis due to
R\"o{\ss}ler~\cite{Roe04d,Roe03}, we derive conditions for the
coefficients of the continuous SRK method assuring weak order one
and two, respectively. As a result of this analysis, we obtain
order conditions for the coefficients of the CSRK
method~(\ref{SRK-method-Ito-Wm-allg01}) which coincide for
$\theta=1$ with the conditions stated in \cite{Roe04e,Roe03}. The
following theorem for continuous SRK methods is an extension of
Theorem~5.1.1 in~\cite{Roe03}.

\begin{The} \label{SRK-theorem-ito-ord2-Wm-main1}
    Let $a^i, b^{ij} \in C_P^{3,6}(I \times
    \mathbb{R}^d,\mathbb{R})$ for
    $1 \leq i \leq d$, $1 \leq j \leq m$.
    If the coefficients of the continuous stochastic Runge--Kutta
    method (\ref{SRK-method-Ito-Wm-allg01})
    fulfill the equations
    \begin{alignat*}{5}
        1&. \,\,\,\, \alpha(\theta)^T e = \theta \quad \quad \quad
        &2. \,\,\,\, &{\beta^{(4)}}(\theta)^T e = 0 \quad \quad \quad
        &3. \,\,\,\, &{\beta^{(3)}}(\theta)^T e = 0 \\
        4&. \,\,\,\, ({\beta^{(1)}}(\theta)^T e)^2 = \theta \quad \quad \quad
        &5. \,\,\,\, &{\beta^{(2)}}(\theta)^T e = 0 \quad \quad \quad
        &6. \,\,\,\, &{\beta^{(1)}}(\theta)^T {B^{(1)}} e = 0 \\
        7&. \,\,\,\, {\beta^{(3)}}(\theta)^T {B^{(2)}} e = 0
    \end{alignat*}
    for $\theta \in [0,1]$ and the equations
    {\allowdisplaybreaks
    \begin{alignat*}{3}
        8&. \quad \alpha(1)^T A^{(0)} e = \tfrac{1}{2}
        \qquad \qquad
        &9. \quad &\alpha(1)^T (B^{(0)} e)^2 = \tfrac{1}{2} \\
        10&. \quad ({\beta^{(1)}}(1)^T e) (\alpha^T B^{(0)} e) =
        \tfrac{1}{2}
        \qquad \qquad
        &11. \quad &({\beta^{(1)}}(1)^T e) ({\beta^{(1)}}(1)^T A^{(1)} e) = \tfrac{1}{2} \\
        12&. \quad {\beta^{(3)}}(1)^T A^{(2)} e = 0
        \qquad \qquad
        &13. \quad &{\beta^{(2)}}(1)^T B^{(1)} e = 1 \\
        14&. \quad {\beta^{(4)}}(1)^T B^{(2)} e = 1
        \qquad \qquad
        &15. \quad &({\beta^{(1)}}(1)^T e) ({\beta^{(1)}}(1)^T
        (B^{(1)} e)^2) = \tfrac{1}{2} \\
        16&. \quad ({\beta^{(1)}}(1)^T e) ({\beta^{(3)}}(1)^T
        (B^{(2)} e)^2) = \tfrac{1}{2} \quad
        &17. \quad &{\beta^{(1)}}(1)^T (B^{(1)} (B^{(1)} e)) =
        0 \\
        18&. \quad {\beta^{(3)}}(1)^T (B^{(2)}
        (B^{(1)}e)) = 0 \qquad
        &19. \quad &{\beta^{(3)}}(1)^T (A^{(2)} (B^{(0)}
        e)) = 0 \\
        20&. \quad {\beta^{(1)}}(1)^T (A^{(1)} (B^{(0)} e)) =
        0 \qquad
        &21. \quad &\alpha(1)^T (B^{(0)} (B^{(1)} e)) = 0 \\
        22&. \quad {\beta^{(2)}}(1)^T A^{(1)} e = 0 \qquad
        &23. \quad &{\beta^{(4)}}(1)^T A^{(2)} e = 0 \\
        24&. \quad {\beta^{(1)}}(1)^T ((A^{(1)} e) (B^{(1)} e)) = 0
        \qquad
        &25. \quad &{\beta^{(3)}}(1)^T ((A^{(2)} e)
        (B^{(2)} e)) = 0 \\
        26&. \quad {\beta^{(4)}}(1)^T (A^{(2)} (B^{(0)}
        e)) = 0 \qquad
        &27. \quad &{\beta^{(2)}}(1)^T (A^{(1)} (B^{(0)} e))
        = 0 \\
        28&. \quad {\beta^{(2)}}(1)^T (A^{(1)} (B^{(0)} e)^2) = 0
        \qquad
        &29. \quad &{\beta^{(4)}}(1)^T (A^{(2)} (B^{(0)}
        e)^2) = 0 \\
        30&. \quad {\beta^{(3)}}(1)^T (B^{(2)} (A^{(1)}
        e)) = 0 \qquad
        &31. \quad &{\beta^{(1)}}(1)^T (B^{(1)} ( A^{(1)} e)) =
        0 \\
        32&. \quad {\beta^{(2)}}(1)^T (B^{(1)} e)^2 = 0 \qquad
        &33. \quad &{\beta^{(4)}}(1)^T (B^{(2)} e)^2 = 0 \\
        34&. \quad {\beta^{(4)}}(1)^T (B^{(2)} (B^{(1)}
        e)) = 0 \qquad
        &35. \quad &{\beta^{(2)}}(1)^T (B^{(1)} (B^{(1)} e))
        = 0 \\
        36&. \quad {\beta^{(1)}}(1)^T (B^{(1)} e)^3 = 0 \qquad
        &37. \quad &{\beta^{(3)}}(1)^T (B^{(2)} e)^3 = 0 \\
        38&. \quad {\beta^{(1)}}(1)^T (B^{(1)} (B^{(1)}
        e)^2) = 0 \qquad
        &39. \quad &{\beta^{(3)}}(1)^T (B^{(2)} (B^{(1)}
        e)^2) = 0
    \end{alignat*}
}
    \begin{alignat*}{3}
        40&. \quad \alpha(1)^T ((B^{(0)} e) (B^{(0)}
        (B^{(1)} e))) = 0 \\
        41&. \quad {\beta^{(1)}}(1)^T ((A^{(1)} (B^{(0)}
        e)) (B^{(1)} e)) = 0 \\
        42&. \quad {\beta^{(3)}}(1)^T ((A^{(2)}
        (B^{(0)} e)) (B^{(2)} e)) = 0 \\
        43&. \quad {\beta^{(1)}}(1)^T (A^{(1)} (B^{(0)} (B^{(1)} e))) = 0 \\
        44&. \quad {\beta^{(3)}}(1)^T (A^{(2)} (B^{(0)}
        (B^{(1)} e))) = 0 \\
        45&. \quad {\beta^{(1)}}(1)^T (B^{(1)} (A^{(1)}
        (B^{(0)} e))) = 0 \\
        46&. \quad {\beta^{(3)}}(1)^T (B^{(2)} (A^{(1)}
        (B^{(0)} e))) = 0 \\
        47&. \quad {\beta^{(1)}}(1)^T ((B^{(1)} e) (B^{(1)}
        (B^{(1)} e))) = 0 \\
        48&. \quad {\beta^{(3)}}(1)^T ((B^{(2)} e)
        (B^{(2)} (B^{(1)} e))) = 0 \\
        49&. \quad {\beta^{(1)}}(1)^T (B^{(1)} (B^{(1)}
        (B^{(1)} e))) = 0 \\
        50&. \quad {\beta^{(3)}}(1)^T (B^{(2)} (B^{(1)}
        (B^{(1)} e))) = 0
    \end{alignat*}
    are fulfilled then the continuous stochastic Runge--Kutta
    method~(\ref{SRK-method-Ito-Wm-allg01}) converges with order
    2 in the weak sense.
\end{The}

{\textbf{Proof.}} This follows directly from the order conditions
for stochastic Runge--Kutta methods in Theorem~5.1.1
in~\cite{Roe03} with the weights replaced by some continuous
functions and taking into account $\theta \, h$ instead of $h$ for
the expansion of the solution. Considering the order 1 conditions
for the time discrete case, we obtain for example $\alpha^T e \, h
= h$ and $({\beta^{(1)}}^T e)^2 \, h = h$. Here, the left hand
side results from the expansion of the approximation process while
the right hand side comes from the expansion of the solution
process at time $t+h$. Now, in the continuous time case, in order
to fulfill condition (\ref{St-lg-local-error-estimation:b}) of
Theorem~\ref{St-lg-theorem-main} we replace the weights of the SRK
method by some continuous functions depending on the parameter
$\theta$ and we consider the solution at time $t+\theta \, h$.
Thus, we obtain the order conditions $\alpha(\theta)^T e \, h =
\theta \, h$ and $(\beta^{(1)}(\theta)^T e)^2 \, h = \theta \, h$.
The remaining order 1 conditions can be obtained in the same
manner. However, the order 2 conditions need to be fulfilled only
at time $t+h$, i.e., for $\theta=1$, due to
(\ref{St-lg-local-error-estimation:a}) of
Theorem~\ref{St-lg-theorem-main}. Thus, we arrive directly at the
conditions of Theorem~\ref{SRK-theorem-ito-ord2-Wm-main1}.
Further, condition (\ref{St-lg-cond2}) of
Theorem~\ref{St-lg-theorem-main} is fulfilled for the
approximations $Y(t)$ at time $t \in I_h$, see \cite{Roe03}.
\hfill $\Box$

\begin{Bem} \label{Bemerkung-1}
    We have to solve 50 equations for $m > 1$ for schemes of order 2.
    However in case of $m=1$ the
    50 conditions are reduced to 28 conditions (see, e.g.
    \cite{Roe04e,Roe03}). Further, explicit CSRK
    methods of order 2 need $s \geq 3$ stages. This is due to the
    conditions 6.\ and 15., which can not be fulfilled for an explicit
    CSRK scheme with $s \leq 2$ stages.
\end{Bem}

\begin{Bem} \label{Bemerkung-2}
    The conditions 1.\--7.\ of
    Theorem~\ref{SRK-theorem-ito-ord2-Wm-main1} for $\theta=1$
    are exactly the order conditions for continuous SRK methods
    of weak order 1.
\end{Bem}

As for time discrete SRK schemes, we distinguish between the
stochastic and the deterministic order of convergence of CSRK
schemes. Let $p_S=p$ denote the order of convergence of the CSRK
method if it is applied to an SDE and let $p_D$ with $p_D \geq
p_S$ denote the order of convergence of the CSRK method if it is
applied to a deterministic ordinary differential equation (ODE),
i.e., SDE~(\ref{St-lg-sde-ito-1}) with $b \equiv 0$. Then, we
write $(p_D,p_S)$ for the orders of convergence in the following
\cite{Roe03}.

Next, we want to calculate coefficients for the CSRK method
(\ref{SRK-method-Ito-Wm-allg01}) which fulfill the order
conditions of Theorem~\ref{SRK-theorem-ito-ord2-Wm-main1}. Since
the conditions of Theorem~\ref{SRK-theorem-ito-ord2-Wm-main1}
coincide for $\theta=1$ with the order conditions of the
underlying discrete time SRK method, we can extend any SRK scheme
to a continuous SRK scheme. To this end, we need some weight
functions depending on $\theta \in [0,1]$ which coincide for
$\theta=1$ with the weights of the discrete time SRK scheme as
boundary conditions. Then, the conditions 1.\--50.\ of
Theorem~\ref{SRK-theorem-ito-ord2-Wm-main1} are fulfilled for
$\theta=1$. Further, the weight functions have to fulfill the
conditions 1.\--7.\ of Theorem~\ref{SRK-theorem-ito-ord2-Wm-main1}
also for all $\theta \in [0,1]$. For $\theta=0$, the right hand
side of (\ref{SRK-method-Ito-Wm-allg01}) has to be equal to
$Y(t_n)$, i.\,e., the weight functions have to vanish for
$\theta=0$, which yields the second boundary conditions. As a
result of this, if coefficients for a discrete time SRK scheme are
already given, then the calculation of coefficients for a CSRK
method reduces to the determination of suitable weight functions
fulfilling some boundary conditions.

As a first example, we consider the order $(1,1)$ SRK scheme with
$s=1$ stage having the weights
$\bar{\alpha}_1=\bar{\beta}_1^{(1)}=1$ and $\bar{\beta}_1^{(2)} =
\bar{\beta}_1^{(3)} = \bar{\beta}_1^{(4)} = 0$. This is the well
known Euler-Maruyama scheme. For the continuous extension, we have
to replace the weights $\bar{\alpha}_1$ and $\bar{\beta}_1^{(k)}$
by some weight functions $\alpha_1, \beta_1^{(k)} \in
C([0,1],\mathbb{R})$ for $k=1,\ldots,4$. Due to
Remark~\ref{Bemerkung-2}, we only have the boundary conditions
$\alpha_1(0)=\beta_1^{(k)}=0$, $\alpha_1(1) = \beta_1^{(1)}(1) =
1$ and $\beta_1^{(2)}(1) = \beta_1^{(3)}(1) = \beta_1^{(4)}(1) =
0$ for $k=1, \ldots, 4$, which have to be fulfilled for an order
(1,1) CSRK scheme. Therefore, we can choose $\alpha_1(\theta)=
\beta_1^{(1)}(\theta) = \theta$ and $\beta_1^{(2)}(\theta) =
\beta_1^{(3)}(\theta) = \beta_1^{(4)}(\theta) = 0$ for $\theta \in
[0,1]$. This results in the linearly interpolated Euler-Maruyama
scheme, which is a CSRK scheme of uniform order $(1,1)$.

Because there are still some degrees of freedom in choosing the
continuous functions $\alpha_1, \beta_1^{(k)} \in
C([0,1],\mathbb{R})$ for $k=1,\ldots,4$, we can try to find some
optimal functions in the sense that additional order conditions
are satisfied. This results in CSRK schemes with usually smaller
error constants in the truncation error. Therefore, these schemes
are called optimal CSRK schemes in the following. For the already
considered example of the order $(1,1)$ CSRK scheme with $s=1$
stage, we can calculate an optimal scheme as follows: firstly,
considering the continuous condition 1.\ of
Theorem~\ref{SRK-theorem-ito-ord2-Wm-main1}, we obtain for $\theta
\in [0,1]$ the additional condition
\begin{equation*}
    \alpha_1(\theta) = \theta \, ,
\end{equation*}
and from the continuous condition 4.\ that
\begin{equation} \label{Sec4:condition-b1}
    \beta_1^{(1)}(\theta) = c_1 \sqrt{\theta}
\end{equation}
for some $c_1 \in \{-1,1\}$. Further, the continuous conditions
2., 3.\ and 5.\ result in
\begin{equation} \label{Sec4:conditions-b2}
    \beta_1^{(2)}(\theta) =0 \, , \qquad \beta_1^{(3)}(\theta) =0 \, ,
    \qquad \beta_1^{(4)}(\theta) =0 \, ,
\end{equation}
for all $\theta \in [0,1]$. Now, all continuous weight functions
are uniquely determined for the optimal order $(1,1)$ CSRK scheme,
which is still a continuous extension of the Euler--Maruyama
scheme.

\begin{table}[tbp]
\renewcommand{\arraystretch}{1.0}
\begin{equation*}
\begin{array}{r|ccc|ccc|ccc}
    0 &&&&&&&&& \\
    \frac{2}{3} & \frac{2}{3} & & & \frac{2}{3} &&&&& \\
    \hline
    0 &&&&&&&&& \\
    0 & 0 & &  & 0 &&&&& \\
    \hline
    0 &&&&&&&&& \\
    0 & 0 &&& 0 &&&&& \\
    \hline
    & \theta-\frac{3}{4} \theta^2 & & \frac{3}{4} \theta^2 &
    \sqrt{\theta} & & 0 & 0 & \quad & 0 \\
    \cline{2-10}
    & & & & 0 & & 0 & 0 & & 0
\end{array}
\end{equation*}
\caption{Coefficients of the optimal CSRK scheme CRDI1WM with
$p_D=2.0$ and $p_S=1.0$.} \label{Table-Coeff-CRDI1}
\end{table}

Now, let us consider the order $(2,1)$ SRK scheme RDI1WM proposed
in \cite{DeRoe06a} with two stages for the deterministic part and
one stage for the stochastic one. The coefficients of RDI1WM for
$A^{(0)}$, $B^{(0)}$, $A^{(1)}$, $B^{(1)}$ and $A^{(2)}$,
$B^{(2)}$ are given in Table~\ref{Table-Coeff-CRDI1} and the
corresponding weights are
\begin{equation*}
    \bar{\alpha} = [\tfrac{1}{4}, \tfrac{3}{4}], \qquad \bar{\beta}^{(1)} =
    [1, 0] \, , \qquad \bar{\beta}^{(2)} = \bar{\beta}^{(3)} = \bar{\beta}^{(4)} =
    [0,0] \, .
\end{equation*}
These coefficients are optimal in the sense that additionally some
higher order conditions are fulfilled. For the continuous
extension, we want to preserve the deterministic order 2 of the
scheme. Therefore, beside the order one conditions 1.\--7.\ of
Theorem~\ref{SRK-theorem-ito-ord2-Wm-main1} for $\theta=1$ we have
to take into account the continuous version of the deterministic
order one condition which coincides with condition 1.\ in
Theorem~\ref{SRK-theorem-ito-ord2-Wm-main1}. Thus, we are looking
for some weight functions $\alpha_1, \alpha_2, \beta_1^{(k)} \in
C([0,1],\mathbb{R})$ fulfilling the boundary conditions
$\alpha_1(0) = \alpha_2(0) = \beta_1^{(k)}(0) = 0$,
\begin{equation*}
    \alpha_1(1) = \frac{1}{4} \, , \qquad \alpha_2(1) =
    \frac{3}{4} \, , \qquad \beta_1^{(1)}(1) =
    1 \, ,
\end{equation*}
and $\beta_1^{(2)}(1) = \beta_1^{(3)}(1) = \beta_1^{(4)}(1) = 0$
for $k=1, \ldots, 4$. Further, due to the continuous version of
the deterministic order one condition, the equation
\begin{equation} \label{Sec4:condition-a1}
    \alpha_1(\theta) + \alpha_2(\theta) = \theta
\end{equation}
has to be fulfilled for all $\theta \in [0,1]$. In order to save
computational effort we set $\beta_1^{(k)}(\theta)=0$ for
$k=2,3,4$ for all $\theta \in [0,1]$. Again, there are still some
degrees of freedom in choosing the continuous functions $\alpha_1,
\alpha_2, \beta_1^{(1)} \in C([0,1],\mathbb{R})$. Now, we can
proceed exactly in the same way as for the order $(1,1)$ CSRK
schemes. Considering additionally the continuous order conditions
2., 3., 4.\ and 5., we obtain again (\ref{Sec4:condition-b1}) and
(\ref{Sec4:conditions-b2}). Finally, we calculate from the
continuous version of condition 8.\ that
\begin{equation*}
    \alpha_2(\theta) = \frac{3}{4} \theta^2
\end{equation*}
and thus $\alpha_1(\theta) = \theta-\tfrac{3}{4} \theta^2$ due to
(\ref{Sec4:condition-a1}) for all $\theta \in [0,1]$. The
coefficients of the CSRK scheme CRDI1WM are given in
Table~\ref{Table-Coeff-CRDI1}.

\begin{table}[tbp]
\renewcommand{\arraystretch}{1.0}
\begin{equation*}
\begin{array}{r|ccccc|ccccc|cccccc}
    0 & & & &&  &&&&&  &&& \\
    1 & 1 & &&&  & 1 & &&  &&& \\
    0 & 0 & & 0 &&  & 0 & & 0 &  &&& \\
    \hline
    0 &&&  &&&&&  &&& \\
    \frac{2}{3} & \frac{2}{3} & & & &  & \sqrt{\frac{2}{3}} &&&&  &&& \\
    \frac{2}{3} & \frac{2}{3} & & 0 & &  & -\sqrt{\frac{2}{3}} && 0 &&  &&& \\
    \hline
    0 & & &&&  &&&&&  &&& \\
    0 & 0 & &&&  & \sqrt{2} &&&&  &&& \\
    0 & 0 & & 0 & & & -\sqrt{2} & & 0 & & &&& \\
    \hline
    & \theta-\frac{1}{2} \theta^2 & & \frac{1}{2} \theta^2 & & 0 &
    \sqrt{\theta} - \frac{3}{4} \theta^{\frac{3}{2}} & & \frac{3}{8} \theta^{\frac{3}{2}}
    & & \frac{3}{8} \theta^{\frac{3}{2}}
    & & 0 & & \frac{\sqrt{6}}{4} \theta & & -\frac{\sqrt{6}}{4} \theta \\
    \cline{2-17}
    & & & & & & -\frac{1}{4} \theta^{\frac{3}{2}} & & \frac{1}{8} \theta^{\frac{3}{2}}
    & & \frac{1}{8} \theta^{\frac{3}{2}}
    & & 0 & & \frac{\sqrt{2}}{4} \theta & & -\frac{\sqrt{2}}{4} \theta
\end{array}
\end{equation*}
\caption{Coefficients of the optimal CSRK scheme CRDI2WM with
$p_D=2.0$ and $p_S=2.0$.} \label{Table-Coeff-CRDI2}
\end{table}

We consider now the continuous extension of the order $(2,2)$ SRK
scheme RDI2WM with $s=3$ stages considered in \cite{DeRoe06a}. The
coefficients of the SRK scheme RDI2WM for $A^{(0)}$, $B^{(0)}$,
$A^{(1)}$, $B^{(1)}$ and $A^{(2)}$, $B^{(2)}$ can be found in
Table~\ref{Table-Coeff-CRDI2} and the weights are given as
\begin{alignat*}{5}
    \bar{\alpha} = [\tfrac{1}{2}, \tfrac{1}{2}, 0], \qquad
    \bar{\beta}^{(1)} &=
    [\tfrac{1}{4}, \tfrac{3}{8}, \tfrac{3}{8}] \, , &\qquad
    \bar{\beta}^{(2)} &= [0, \tfrac{\sqrt{6}}{4}, -\tfrac{\sqrt{6}}{4}] \,
    , \\
    \qquad \bar{\beta}^{(3)} &= [-\tfrac{1}{4}, \tfrac{1}{8}, \tfrac{1}{8}] \, ,
    &\qquad \bar{\beta}^{(4)} &= [0, \tfrac{\sqrt{2}}{4}, -\tfrac{\sqrt{2}}{4}] \, .
\end{alignat*}
We proceed in the same way as for the scheme CRDI1WM by taking
into account some additional order conditions. As the stage number
of the deterministic part of the scheme is only two, we set
$\alpha_3(\theta) = 0$.
From the order 1 conditions of
Theorem~\ref{SRK-theorem-ito-ord2-Wm-main1} follows for $\alpha_i,
\beta_i^{(k)} \in C([0,1],\mathbb{R})$, $1 \leq i \leq 3$, $1 \leq
k \leq 4$, that
\begin{alignat*}{5}
    \alpha_1(\theta) &= \theta -\alpha_2(\theta) ,
    &\quad \quad \quad \quad & \\
    \beta^{(1)}_1(\theta) &= c_1 \sqrt{\theta} -2 \beta^{(1)}_2(\theta) , &\quad
    \beta^{(1)}_2(\theta) &= \beta^{(1)}_3(\theta) , \\
    \beta^{(3)}_1(\theta) &= -2 \beta^{(3)}_2(\theta) , &\quad
    \beta^{(3)}_2(\theta) &= \beta^{(3)}_3(\theta) , \\
    \beta^{(2)}_1(\theta) &= -\beta^{(2)}_2(\theta) - \beta^{(2)}_3(\theta) , &\quad
    \beta^{(4)}_1(\theta) &= -\beta^{(4)}_2(\theta)-\beta^{(4)}_3(\theta) .
\end{alignat*}
Further, the functions $\alpha_i$ and $\beta_i^{(k)}$ have to
fulfill for $\theta=0$ and $\theta=1$ the boundary conditions
\begin{alignat*}{5}
    \alpha_2(1) &= \frac{1}{2}, \qquad &\beta_2^{(1)}(1) &=
    \frac{3}{8}, \qquad &\beta_2^{(2)}(1) &= \frac{\sqrt{6}}{4}, \qquad
    &\beta_3^{(2)}(1) &= -\frac{\sqrt{6}}{4}, \\
    \beta_2^{(3)}(1) &= \frac{1}{8} , \qquad &\beta_2^{(4)}(1) &=
    \frac{\sqrt{2}}{4}, \qquad &\beta_3^{(4)}(1) &= -\frac{\sqrt{2}}{4} & &
\end{alignat*}
and $\alpha_i(0) = \beta_i^{(k)}(0) = 0$ for $1 \leq i \leq 3$ and
$1 \leq k \leq 4$.
Here, the continuous versions of the order conditions 8.--50.\ of
Theorem~\ref{SRK-theorem-ito-ord2-Wm-main1}, i.e., conditions
8.--50.\ with $\alpha(\theta)$ and $\beta^{(k)}(\theta)$ instead
of $\alpha(1)$ and $\beta^{(k)}(1)$ for $k=1, \ldots, 4$, and with
appropriate right hand sides are automatically fulfilled except of
conditions 8., 9., 10., 11., 13., 14., 15., 16., 22., 32.\ and 33.
The continuous versions of condition 8.\ and 9.\ yield
\begin{equation} \label{cond_08}
    \alpha_2(\theta) = \frac{1}{2} \, \theta^2 .
\end{equation}
However, condition 10.\ alternatively yields
\begin{equation*}
    \alpha_2(\theta) = \frac{1}{2} \, \theta^{\frac{3}{2}},
\end{equation*}
so one has to decide whether 8.\ and 9.\ or 10.\ should be
fulfilled. The conditions 22.\ and 32.\ coincide and provide
\begin{equation} \label{cond_01}
    \beta_3^{(2)}(\theta) = -\beta_2^{(2)}(\theta) .
\end{equation}
Further, condition 33.\ results in
\begin{equation} \label{cond_02}
    \beta_3^{(4)}(\theta) = - \beta_2^{(4)}(\theta) .
\end{equation}
Taking into account condition 13., we calculate that
\begin{equation} \label{cond_03}
    \beta_2^{(2)}(\theta) = \beta_3^{(2)}(\theta) + \sqrt{\frac{3}{2}} \, \theta ,
\end{equation}
and due to condition 14.\ that
\begin{equation} \label{cond_04}
    \beta_2^{(4)}(\theta) = \beta_3^{(4)}(\theta) + \frac{1}{\sqrt{2}} \, \theta.
\end{equation}
Now, if we combine conditions 32.\ and 33.\ with 13.\ and 14.\
then we can determine $\beta_3^{(2)}$ and $\beta_3^{(4)}$ uniquely
as
\begin{equation} \label{cond_05}
    \beta_3^{(2)}(\theta) = -\frac{\sqrt{3}}{2 \sqrt{2}} \, \theta, \qquad \qquad
    \beta_3^{(4)}(\theta) = -\frac{1}{2 \sqrt{2}} \, \theta .
\end{equation}
Considering condition 11., which coincides with 15., one obtains
\begin{equation} \label{cond_06}
    \beta_2^{(1)}(\theta) = \frac{3}{8} \, \theta^{\frac{3}{2}} .
\end{equation}
Finally, condition 16.\ yields that
\begin{equation} \label{cond_07}
    \beta_2^{(3)}(\theta) = \frac{1}{8} \, \theta^{\frac{3}{2}} .
\end{equation}
Thus, one can choose from these additional conditions in order to
minimize the truncation constant. Especially, if we consider
equations (\ref{cond_08})--(\ref{cond_07}) then we obtain the
scheme CRDI2WM presented in Table~\ref{Table-Coeff-CRDI2}.

% -----------------------------------------------------
%
\begin{table}[tbp]
\renewcommand{\arraystretch}{1.0}
\begin{equation*}
\begin{array}{r|ccccc|ccccc|cccccc}
    0 & & & &&  &&&&&  &&& \\
    \frac{1}{2} & \frac{1}{2} & &&&  & \frac{9-2\sqrt{15}}{14} & &&  &&& \\
    \frac{3}{4} & 0 & & \frac{3}{4} &&  & \frac{18+3 \sqrt{15}}{28} & & 0 &  &&& \\
    \hline
    0 &&&  &&&&&  &&& \\
    \frac{2}{3} & \frac{2}{3} & & & &  & \sqrt{\frac{2}{3}} &&&&  &&& \\
    \frac{2}{3} & \frac{2}{3} & & 0 & &  & -\sqrt{\frac{2}{3}} && 0 &&  &&& \\
    \hline
    0 & & &&&  &&&&&  &&& \\
    0 & 0 & &&&  & \sqrt{2} &&&&  &&& \\
    0 & 0 & & 0 & & & -\sqrt{2} & & 0 & & &&& \\
    \hline
    & \theta-\frac{7}{9} \theta^2 & & \frac{1}{3} \theta^2 & & \frac{4}{9} \theta^2 &
    \sqrt{\theta} - \frac{3}{4} \theta^{\frac{3}{2}} & & \frac{3}{8} \theta^{\frac{3}{2}}
    & & \frac{3}{8} \theta^{\frac{3}{2}}
    & & 0 & & \frac{\sqrt{6}}{4} \theta & & -\frac{\sqrt{6}}{4} \theta \\
    \cline{2-17}
    & & & & & & -\frac{1}{4} \theta^{\frac{3}{2}} & & \frac{1}{8} \theta^{\frac{3}{2}}
    & & \frac{1}{8} \theta^{\frac{3}{2}}
    & & 0 & & \frac{\sqrt{2}}{4} \theta & & -\frac{\sqrt{2}}{4} \theta
\end{array}
\end{equation*}
\caption{Coefficients of the optimal CSRK scheme CRDI3WM with
$p_D=3$ and $p_S=2$.} \label{Table-Coeff-CRDI3}
\end{table}
If we allow now three stages for the deterministic part, i.e.,
$\alpha_3 \neq 0$, then we can construct explicit CSRK schemes of
order $p_D=3$ and $p_S=2$. Therefore, we extend the order $(3,2)$
explicit SRK scheme RDI3WM introduced in \cite{DeRoe06a} with
coefficients $A^{(0)}$, $B^{(0)}$, $A^{(1)}$, $B^{(1)}$ and
$A^{(2)}$, $B^{(2)}$ from Table~\ref{Table-Coeff-CRDI3} and with
the weights
\begin{alignat*}{5}
    \bar{\alpha} = [\tfrac{2}{9}, \tfrac{1}{3}, \tfrac{4}{9}], \qquad
    \bar{\beta}^{(1)} &=
    [\tfrac{1}{4}, \tfrac{3}{8}, \tfrac{3}{8}] \, , &\qquad
    \bar{\beta}^{(2)} &= [0, \tfrac{\sqrt{6}}{4}, -\tfrac{\sqrt{6}}{4}] \,
    , \\
    \qquad \bar{\beta}^{(3)} &= [-\tfrac{1}{4}, \tfrac{1}{8}, \tfrac{1}{8}] \, ,
    &\qquad \bar{\beta}^{(4)} &= [0, \tfrac{\sqrt{2}}{4}, -\tfrac{\sqrt{2}}{4}] \, .
\end{alignat*}
Here, the well known order three conditions
\begin{equation*}
    \bar{\alpha}^T (A^{(0)} e)^2 = \frac{1}{3} \, , \quad \quad \quad \quad
    \bar{\alpha}^T (A^{(0)} (A^{(0)} e)) = \frac{1}{6} \, ,
\end{equation*}
for deterministic ODEs are fulfilled \cite{HNW93}. Obviously, the
scheme RDI3WM differs from RDI2WM only in $A^{(0)}$, $B^{(0)}$ and
$\bar{\alpha}$. Therefore, we choose the functions $\beta_i^{(k)}
\in C([0,1],\mathbb{R})$ for $1 \leq i \leq 3$ and $1 \leq k \leq
4$ in the same way as for the CSRK scheme CRDI2WM.
Thus, we only have to specify the functions $\alpha_i \in
C([0,1],\mathbb{R})$ for $1 \leq i \leq 3$. For order $p_D=3$, the
continuous version of condition 8.\ in
Theorem~\ref{SRK-theorem-ito-ord2-Wm-main1}
\begin{equation*} \label{continuous-vers-cond-6}
    \alpha(\theta)^T A^{(0)} e = \frac{1}{2} \theta^2
\end{equation*}
has to be fulfilled as well. As a result of this, we yield the
conditions
\begin{alignat*}{5}
    \alpha_1(\theta) &= \frac{1}{2} \alpha_3(\theta) +
    \theta - \theta^2, &\qquad
    \alpha_2(\theta) &= \theta^2 - \frac{3}{2} \alpha_3(\theta),
\end{alignat*}
with boundary conditions $\alpha_3(0)=0$ and $\alpha_3(1) =
\tfrac{4}{9}$.

In order to specify $\alpha_3 \in C([0,1],\mathbb{R})$, we
consider the continuous version of condition 9.\ which has the
unique solution
\begin{equation} \label{Cond-stetig-9a}
    \alpha_3(\theta) = \alpha_3(1) \, \theta^2.
\end{equation}
Alternatively, we can consider the continuous version of condition
10.\ which yields
\begin{equation*}
    \alpha_3(\theta) =
    \frac{2}{9} \frac{\theta^2(2\sqrt{15}-9)+7\theta^{\frac32}}{\sqrt{15}-1}
\end{equation*}

However, instead of considering condition 9. or 10., we can also
take into account the continuous version of the deterministic
order 3 condition
\begin{equation*}
    \alpha(\theta)^T (A^{(0)} (A^{(0)} e)) =
    \frac{1}{6}
\end{equation*}
which yields
\begin{equation*}
    \alpha_3(\theta) = \alpha_3(1) \, \theta^3 \, .
\end{equation*}

On the other hand, if we consider the continuous version of the
deterministic order 3 condition
\begin{equation*} %\label{condition-glz1}
    \alpha(\theta)^T (A^{(0)} e)^2 = \frac{1}{3}
    \, \theta^3
\end{equation*}
then we get for the deterministic part of the continuous SRK
scheme that
\begin{equation} \label{condition-glz1-eqn}
    \alpha_3(\theta) = \frac{16}{9} \theta^3 - \frac{4}{3}
    \theta^2 \, .
\end{equation}
So, one can choose from these additional conditions in order to
minimize the truncation error of the considered scheme. For the
continuous extension of the SRK scheme RDI3WM, we choose condition
(\ref{Cond-stetig-9a}) which yields the coefficients of the CSRK
scheme CRDI3WM presented in Table~\ref{Table-Coeff-CRDI3}.

% ---------------------------------------------------
%
\begin{table}[tbp]
\renewcommand{\arraystretch}{1.0}
\begin{equation*}
\begin{array}{r|ccccc|ccccc|cccccc}
    0 & & & &&  &&&&&  &&& \\
    \frac{1}{2} & \frac{1}{2} & &&&  & \frac{6-\sqrt{6}}{10} & &&  &&& \\
    1 & -1 & 2 &  &&  & \frac{3+2\sqrt{6}}{5} & & 0 &  &&& \\
    \hline
    0 &&&  &&&&&  &&& \\
    \frac{2}{3} & \frac{2}{3} & & & &  & \sqrt{\frac{2}{3}} &&&&  &&& \\
    \frac{2}{3} & \frac{2}{3} & & 0 & &  & -\sqrt{\frac{2}{3}} && 0 &&  &&& \\
    \hline
    0 & & &&&  &&&&&  &&& \\
    0 & 0 & &&&  & \sqrt{2} &&&&  &&& \\
    0 & 0 & & 0 & & & -\sqrt{2} & & 0 & & &&& \\
    \hline
    & \theta-\frac{5}{6} \theta^2 & & \frac{2}{3} \theta^2 & & \frac{1}{6} \theta^2 &
    \sqrt{\theta} - \frac{3}{4} \theta^{\frac{3}{2}} & & \frac{3}{8} \theta^{\frac{3}{2}}
    & & \frac{3}{8} \theta^{\frac{3}{2}}
    & & 0 & & \frac{\sqrt{6}}{4} \theta & & -\frac{\sqrt{6}}{4} \theta \\
    \cline{2-17}
    & & & & & & -\frac{1}{4} \theta^{\frac{3}{2}} & & \frac{1}{8} \theta^{\frac{3}{2}}
    & & \frac{1}{8} \theta^{\frac{3}{2}}
    & & 0 & & \frac{\sqrt{2}}{4} \theta & & -\frac{\sqrt{2}}{4} \theta
\end{array}
\end{equation*}
\caption{Coefficients of the optimal CSRK scheme CRDI4WM with
$p_D=3$ and $p_S=2$.} \label{Table-Coeff-CRDI4}
\end{table}
Analogously to the procedure for the CSRK scheme CRDI3WM, we
extend the order $(3,2)$ explicit SRK scheme RDI4WM
\cite{DeRoe06a} with coefficients $A^{(0)}$, $B^{(0)}$, $A^{(1)}$,
$B^{(1)}$ and $A^{(2)}$, $B^{(2)}$ from
Table~\ref{Table-Coeff-CRDI4} and with the weights
\begin{alignat*}{5}
    \bar{\alpha} = [\tfrac{1}{6}, \tfrac{2}{3}, \tfrac{1}{6}], \qquad
    \bar{\beta}^{(1)} &=
    [\tfrac{1}{4}, \tfrac{3}{8}, \tfrac{3}{8}] \, , &\qquad
    \bar{\beta}^{(2)} &= [0, \tfrac{\sqrt{6}}{4}, -\tfrac{\sqrt{6}}{4}] \,
    , \\
    \qquad \bar{\beta}^{(3)} &= [-\tfrac{1}{4}, \tfrac{1}{8}, \tfrac{1}{8}] \, ,
    &\qquad \bar{\beta}^{(4)} &= [0, \tfrac{\sqrt{2}}{4}, -\tfrac{\sqrt{2}}{4}] \, .
\end{alignat*}
Then, we obtain the conditions
\begin{alignat*}{5}
    \alpha_1(\theta) &= \alpha_3(\theta) + \theta - \theta^2, &\qquad
    \alpha_2(\theta) &= \theta^2 - 2 \alpha_3(\theta),
\end{alignat*}
with boundary conditions $\alpha_3(0)=0$ and $\alpha_3(1) =
\tfrac{1}{6}$. Choosing again condition (\ref{Cond-stetig-9a})
yields the coefficients of the CSRK scheme CRDI4WM given in
Table~\ref{Table-Coeff-CRDI4}.

However, if we consider condition (\ref{condition-glz1-eqn})
instead of (\ref{Cond-stetig-9a}), then we obtain coefficients for
the CSRK scheme CRDI5WM. The coefficients of the scheme CRDI5WM
coincide with the coefficients of the scheme CRDI4WM, with the
exception of the weights $\alpha_i(\theta)$. From condition
(\ref{condition-glz1-eqn}), we calculate the weights
$\alpha_1(\theta) = \frac{2}{3} \theta^3-\frac{3}{2}
\theta^2+\theta$, $\alpha_2(\theta) = 2 \theta^2-\frac{4}{3}
\theta^3$ and $\alpha_3(\theta) = \frac{2}{3} \theta^3-\frac{1}{2}
\theta^2$ for the scheme CRDI5WM.
\section{Numerical example} \label{Sec:Numerical-Example}
Now, the proposed CSRK scheme CRDI3WM is applied in order to
analyze the empirical order of convergence. The first considered
test equation is a linear SDE (see, e.g., \cite{KP99})
\begin{equation} \label{Simu:linear-SDE1}
    {\mathrm{d}}X(t) = a X(t) \, {\mathrm{d}}t + b X(t) \, {\mathrm{d}}W(t), \qquad X(0)=x_0,
\end{equation}
with $a, b, x_0 \in \mathbb{R}$. In the following, we choose $f(x) =
x$ and we consider the interval $I=[0,2]$. Then the expectation of
the solution at time $t \in I$ is given by $E(X(t)) = x_0 \cdot
\exp(a \, t)$. We consider the case $a=1.5$ with
$b=0.1$ and $x_0=0.1$. \\ \\
%
% ------------------------------------------------------------------------
%
As a second example, a multi-dimensional SDE with a
$2$-dimensional dri\-ving Wiener process is considered:
\begin{equation} \label{Simu:dm-SDE}
    \begin{split}
    {\mathrm{d}} \begin{pmatrix} X^1(t) \\ X^2(t) \end{pmatrix}
    &= \begin{pmatrix} -\frac{273}{512} & 0 \\
    -\frac{1}{160} & -\frac{785}{512} \end{pmatrix} \,
    \begin{pmatrix} X^1(t) \\ X^2(t) \end{pmatrix} \, {\mathrm{d}}t \\
    &+
    \begin{pmatrix} \frac{1}{16} X^1(t) & \frac{1}{16} X^1(t) \\
    \frac{1-2\sqrt{2}}{4} X^2(t) & \frac{1}{10} X^1(t) + \frac{1}{16} X^2(t)
    \end{pmatrix} \, {\mathrm{d}} \begin{pmatrix} W^1(t) \\ W^2(t)
    \end{pmatrix}
    \end{split}
\end{equation}
with initial value $X(0)=(1,1)^T$ and $I=[0,4]$. This SDE system
is of special interest due to the fact that it has non-commutative
noise. Here, we are interested in the second moments which depend
on both, the drift and the diffusion function (see \cite{KP99} for
details). Therefore, we choose $f(x) = (x^1)^2$ and obtain
\begin{equation*}
    E(f(X(t))) = \exp(- t) \, .
\end{equation*}

\begin{figure}[tbp]
\begin{center}
\includegraphics[width=6.8cm]{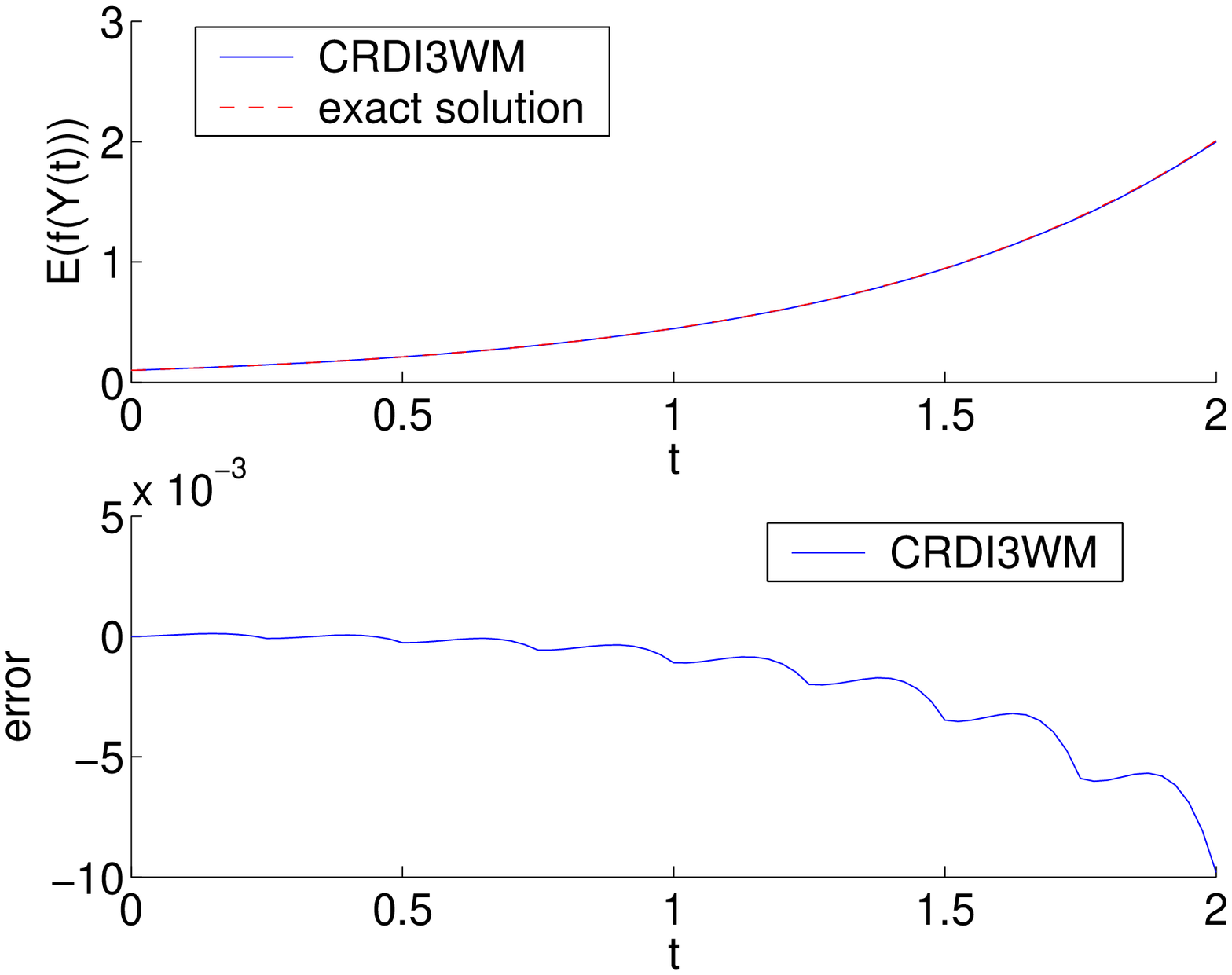}
\includegraphics[width=6.8cm]{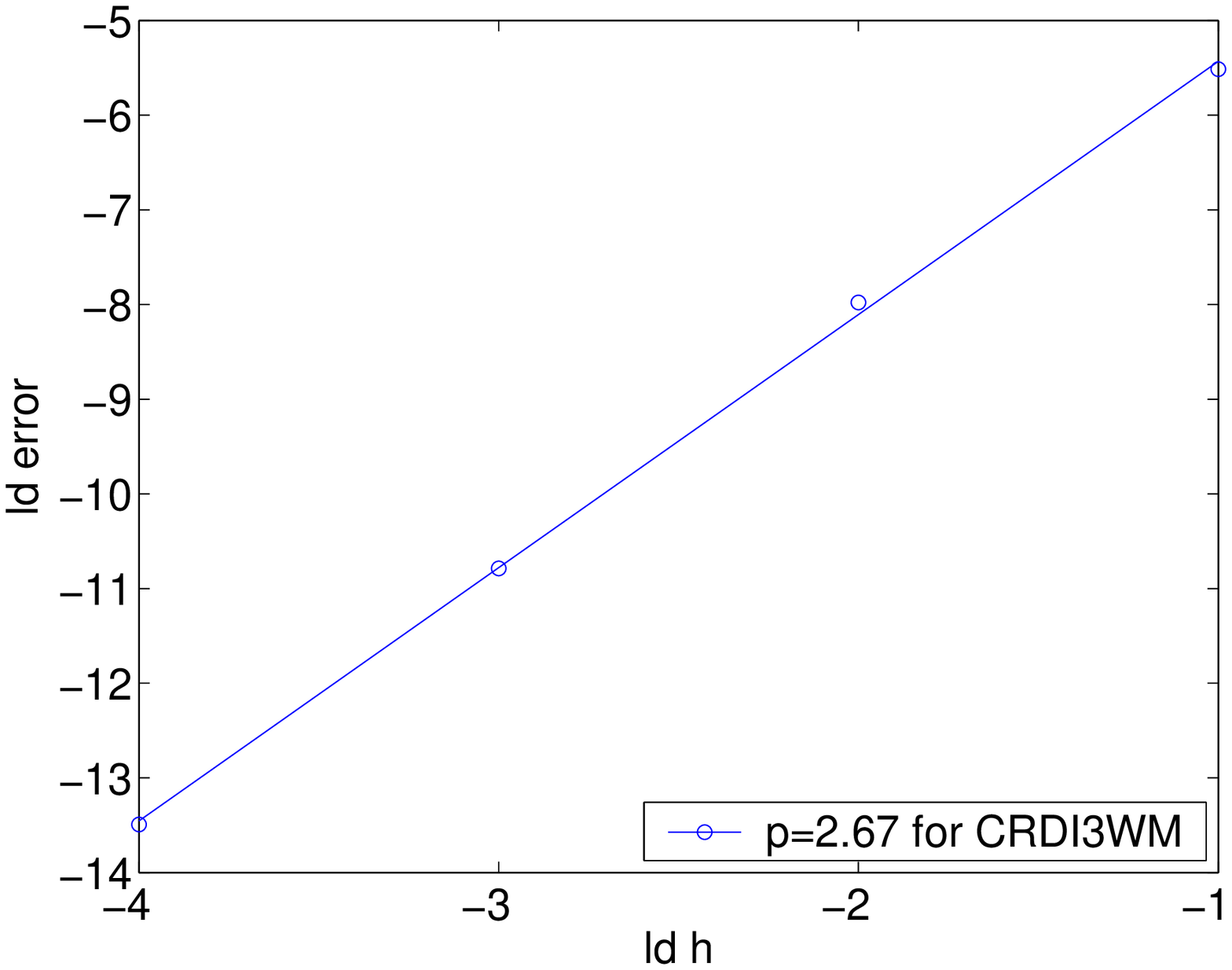}
\caption{Scheme CRDI3WM with SDE~(\ref{Simu:linear-SDE1}).}
\label{Bild2}
\end{center}
\end{figure}

For $t \in I$, we approximate the functional $u = E(f(Y(t)))$ by
Monte Carlo simu\-la\-tion using the sample average $u_{M,h} =
\frac{1}{M} \sum_{m=1}^M f(Y^{(m)}(t))$ of independent simulated
realizations $Y^{(m)}$, $m=1, \ldots, M$, of the considered
approxi\-mation $Y$ and we choose $M=10^9$. Then, the mean error
is given as $\hat{\mu} = u_{M,h} - E(f(X(t)))$ and the estimation
for the variance of the mean error is denoted by
$\hat{\sigma}^2_{\mu}$. Further, a confidence interval $[\hat{\mu}
- \Delta \hat{\mu}, \hat{\mu} + \Delta \hat{\mu} ]$ to the level
of confidence 90\% for the mean error $\mu$ is
calculated~\cite{KP99}. \\ \\
First, the solution $E(f(X(t)))$ is considered as a mapping from
$I$ to $\mathbb{R}$ with $t \mapsto E(f(X(t)))$. Here, the whole
trajectory of the expectation even between the discretization
points has to be determined. Therefore, we apply the CSRK scheme
CRDI3WM with step size $h=0.25$ and determine $E(f(Y({t_n})))$ for
the discretization times $t_n = n \cdot h$ for $n=0,1, \ldots, N$
and the approximation $Y(t)$ is exploited between each pair of
discretization points $t_n$ and $t_{n+1}$ by choosing $\theta \in
\,]0,1[\,$. This has been done for $\theta=0.1, \ldots, 0.9$. The
results are plotted in the left hand side of Figure~\ref{Bild2}
and Figure~\ref{Bild1}. Further, the errors are plotted along the
whole time interval $I$ in the figure below. Here, it turns out
that the continuous extension of the SRK scheme works very well.
%
% ------------------------------------------------------------
%
\setlength{\extrarowheight}{0pt}
\begin{figure}[tbp]
\begin{center}
\includegraphics[width=6.8cm]{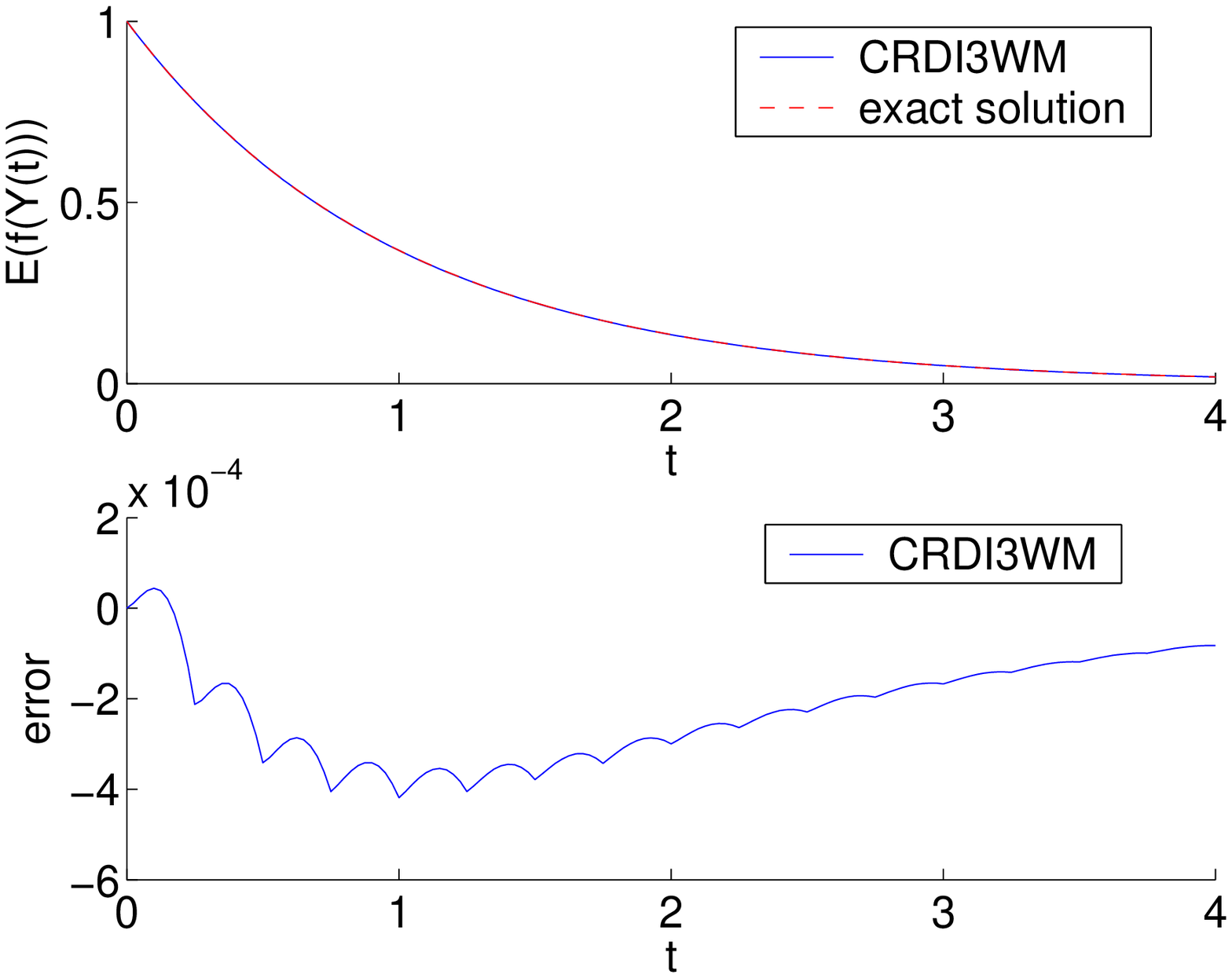}
\includegraphics[width=6.8cm]{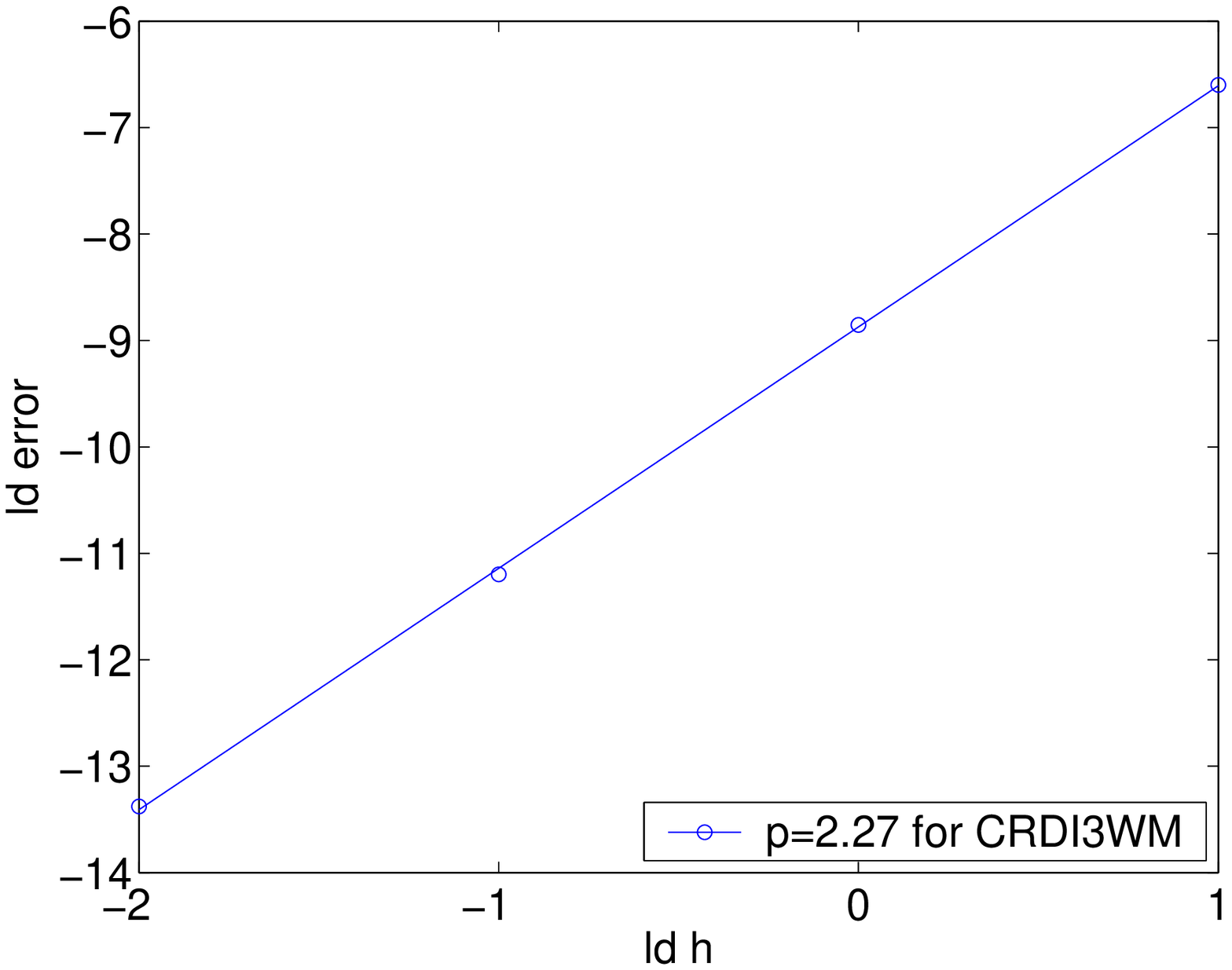}
\caption{Scheme CRDI3WM with SDE~(\ref{Simu:dm-SDE}).}
\label{Bild1}
\end{center}
\end{figure}
\\ \\
Next, SDE~(\ref{Simu:linear-SDE1}) and SDE~(\ref{Simu:dm-SDE}) are
applied for the investigation of the order of convergence.
Therefore, the trajectories are simulated with step sizes $2^{-1},
\ldots, 2^{-4}$ for SDE~(\ref{Simu:linear-SDE1}) and with step
sizes $2^{1}, \ldots, 2^{-2}$ for SDE~(\ref{Simu:dm-SDE}). As an
example, we consider the error $\hat{\mu}$ at time $t=1.7$ for
SDE~(\ref{Simu:linear-SDE1}) and at $t=3.8$ for
SDE~(\ref{Simu:dm-SDE}), which are not discretization points. The
results are plotted on the right hand side of Figure~\ref{Bild2}
and Figure~\ref{Bild1} with double logarithmic scale w.r.t.\ base
two. On the axis of abscissae, the step sizes are plotted against
the errors on the axis of ordinates. Consequently one obtains the
empirical order of convergence as the slope of the printed lines.
In the case of SDE~(\ref{Simu:linear-SDE1}) we get the order $p
\approx 2.67$ and in the case of SDE~(\ref{Simu:dm-SDE}) we get
the order $p \approx 2.27$. Table~\ref{Table1} and
Table~\ref{Table2} contain the corresponding values of the errors,
the variances and the confidence intervals to the level 90\%. The
same has been done for all time points in order to determine the
empirical order of convergence on the whole time interval $I$.
These results are listed in Table~\ref{Table3} for both considered
examples. The very good empirical orders of convergence confirm
our theoretical results for the CSRK scheme CRDI3WM.
\begin{table}
\caption{Scheme CRDI3WM with SDE~(\ref{Simu:linear-SDE1}) at
$t=1.7$.} \label{Table1}
\begin{center}
\renewcommand{\arraystretch}{1.0}
\begin{tabular}{|c|c|c|c|c|}
    \hline
    $h$ & $\hat{\mu}$ & $\hat{\sigma}_{\mu}^2$ & $\hat{\mu} - \Delta \hat{\mu}$ & $\hat{\mu} + \Delta \hat{\mu}$ \\
    \hline
    $2^{-1}$ & -2.188E-02 & 1.146E-09 & -2.189E-02 & -2.187E-02 \\
    $2^{-2}$ & -3.965E-03 & 1.705E-09 & -3.975E-03 & -3.956E-03 \\
    $2^{-3}$ & -5.662E-04 & 1.715E-09 & -5.760E-04 & -5.564E-04 \\
    $2^{-4}$ & -8.682E-05 & 1.746E-09 & -9.672E-05 & -7.691E-05 \\
    \hline
\end{tabular}
\end{center}
\end{table}
%
% ------------------------------------------------------------------------
%
\begin{table}%[htbp]
\caption{Scheme CRDI3WM with SDE~(\ref{Simu:dm-SDE}) at $t=3.8$.}
\label{Table2}
\begin{center}
\renewcommand{\arraystretch}{1.0}
\begin{tabular}{|c|c|c|c|c|}
    \hline
    $h$ & $\hat{\mu}$ & $\hat{\sigma}_{\mu}^2$ & $\hat{\mu} - \Delta \hat{\mu}$ & $\hat{\mu} + \Delta \hat{\mu}$ \\
    \hline
    $2^{1}$ & -1.031E-02 & 9.365E-12 & -1.031E-02 & -1.031E-02 \\
    $2^{0}$ & -2.161E-03 & 2.294E-11 & -2.162E-03 & -2.160E-03 \\
    $2^{-1}$ & -4.258E-04 & 4.453E-11 & -4.274E-04 & -4.242E-04 \\
    $2^{-2}$ & -9.392E-05 & 3.375E-11 & -9.530E-05 & -9.255E-05 \\
    \hline
\end{tabular}
\end{center}
\end{table}
\\
%
% ------------------------------------------------------------------------
%
\begin{table}%[htbp]
\caption{Orders of convergence for CRDI3WM.} \label{Table3}
\begin{center}
\renewcommand{\arraystretch}{1.0}
\begin{tabular}{|c|c||c|c|}
    \hline
    \multicolumn{4}{|c|}{$SDE~(\ref{Simu:linear-SDE1})$} \\
    \hline
    \hline
    $t$ & order & $t$ & order \\
    \hline
    0.1 & 2.35 & 1.1 & 2.94 \\
    0.2 & 3.39 & 1.2 & 2.60 \\
    0.3 & 3.54 & 1.3 & 2.69 \\
    0.4 & 3.16 & 1.4 & 2.67 \\
    0.5 & 2.87 & 1.5 & 2.80 \\
    0.6 & 3.50 & 1.6 & 2.89 \\
    0.7 & 2.32 & 1.7 & 2.67 \\
    0.8 & 2.44 & 1.8 & 2.72 \\
    0.9 & 2.53 & 1.9 & 2.70 \\
    1.0 & 2.79 & 2.0 & 2.80 \\
    \hline
\end{tabular}
\qquad \qquad
\begin{tabular}{|c|c||c|c|}
    \hline
    \multicolumn{4}{|c|}{$SDE~(\ref{Simu:dm-SDE})$} \\
    \hline
    \hline
    $t$ & order & $t$ & order \\
    \hline
    0.2 & 1.82 & 2.2 & 2.30 \\
    0.4 & 1.84 & 2.4 & 2.27 \\
    0.6 & 1.95 & 2.6 & 2.24 \\
    0.8 & 2.07 & 2.8 & 2.24 \\
    1.0 & 2.13 & 3.0 & 2.22 \\
    1.2 & 2.26 & 3.2 & 2.25 \\
    1.4 & 2.32 & 3.4 & 2.26 \\
    1.6 & 2.36 & 3.6 & 2.27 \\
    1.8 & 2.37 & 3.8 & 2.27 \\
    2.0 & 2.31 & 4.0 & 2.24 \\
    \hline
\end{tabular}
\end{center}
\end{table}
\renewcommand{\bf}{}
\renewcommand{\em}{}
\renewcommand{\sc }{}

\end{document}